\documentclass[a4paper,12pt,final]{amsart}
\usepackage{times,a4wide,mathrsfs,amssymb,amsmath,amsthm,enumerate,xypic,tikzsymbols,dsfont}

\newcommand{\C}{\mathbb{C}}

\newcommand{\QQ}{\mathbb{Q}}
\newcommand{\NN}{\mathbb{N}}
\newcommand{\PP}{\mathbb{P}}

\newcommand{\OO}{\mathcal O}

\newcommand{\Sy}{\mathfrak S}

\newcommand{\XX}{\mathcal X}
\newcommand{\YY}{\mathcal Y}

\newcommand{\CC}{\mathcal C}

\newcommand{\EE}{\mathcal E}
\newcommand{\MM}{\mathcal M}

\newcommand{\HH}{\mathcal H}
\newcommand{\FF}{\mathcal F}

\newcommand{\pic}{\hbox{Pic}}

\newcommand{\rom}{\romannumeral}

\newcommand{\one}{\mathds{1}}

\DeclareMathOperator{\ima}{Im}

\DeclareMathOperator{\rank}{rank}
\DeclareMathOperator{\sym}{Sym}
\DeclareMathOperator{\Gr}{Gr}

\newtheorem{theorem}{Theorem}[section]

\newtheorem{lemma}[theorem]{Lemma}
\newtheorem{sublemma}[theorem]{Sublemma}
\newtheorem{corollary}[theorem]{Corollary}
\newtheorem{proposition}[theorem]{Proposition}

\newtheorem{convention}{Conventions}

\newtheorem{nonumbering}{Theorem}

\newtheorem{nonumberingc}{Corollary}

\theoremstyle{definition}
\newtheorem{remark}[theorem]{Remark}
\newtheorem{definition}[theorem]{Definition}

\newtheorem{notation}[theorem]{Notation}

\newtheorem{nonumberingt}{Acknowledgments}

\begin{document}

\author[Robert Laterveer]
{Robert Laterveer}

\address{Institut de Recherche Math\'ematique Avanc\'ee,
CNRS -- Universit\'e 
de Strasbourg,\
7 Rue Ren\'e Des\-car\-tes, 67084 Strasbourg CEDEX,
FRANCE.}
\email{robert.laterveer@math.unistra.fr}

\title{Algebraic cycles and intersections of 2 quadrics}

\begin{abstract} A smooth intersection $Y$ of two quadrics in $\PP^{2g+1}$ has Hodge level 1. We show that such varieties $Y$ have a multiplicative Chow--K\"unneth decomposition, in the sense of Shen--Vial. As a 
consequence, a certain tautological subring of the Chow ring of powers of $Y$ injects into cohomology.
 \end{abstract}

\thanks{Supported by ANR grant ANR-20-CE40-0023.}


\maketitle

\section{Introduction}

Given a smooth projective variety $Y$ over $\C$, let $A^i(Y):=CH^i(Y)_{\QQ}$ denote the Chow groups of $Y$ (i.e. the groups of codimension $i$ algebraic cycles on $Y$ with $\QQ$-coefficients, modulo rational equivalence). The intersection product defines a ring structure on $A^\ast(Y)=\bigoplus_i A^i(Y)$, the {\em Chow ring\/} of $Y$.

Motivated by the particular behaviour of Chow rings of K3 surfaces \cite{BV} and of abelian varieties \cite{Beau}, Beauville \cite{Beau3} has conjectured that for certain special varieties, the Chow ring should admit a multiplicative splitting. To make concrete sense of Beauville's elusive ``splitting property'' conjecture, Shen--Vial \cite{SV} introduced the concept of {\em multiplicative Chow--K\"unneth decomposition\/}. In short, this is a graded decomposition of the Chow motive of a smooth projective variety, such that the intersection product respects the grading (cf. subsection \ref{ss:mck} for details).

It is something of a challenge to understand the class of special varieties admitting a multiplicative Chow--K\"unneth decomposition: abelian varieties, K3 surfaces and cubic hypersurfaces are in this class, hyperelliptic curves are in this class (but not all curves), some Fano varieties are in this class (but not all Fano varieties).
The main result of the present paper aims to contribute to this program:

\begin{nonumbering}[=Theorem \ref{main}] Let $Y\subset\PP^{2g+1}$ be a smooth dimensionally transverse intersection of 2 quadrics. Then $Y$ has a multiplicative Chow--K\"unneth decomposition.
\end{nonumbering}

(We observe that when $Y\subset\PP^{2g}$ is an {\em even-dimensional\/} intersection of 2 quadrics, $Y$ has trivial Chow groups and so in this case the theorem 
is trivially true.)
%

To prove Theorem \ref{main}, we use the {\em Franchetta property\/} for certain universal families, allowing to lift homological relations between generically defined cycles to the level of rational equivalence. To establish the required instances of the Franchetta property, we exploit the connection (via an Abel--Jacobi isomorphism) between $Y$ and the Fano variety of $g-1$-dimensional linear subspaces contained in $Y$. This is similar to (and inspired by) the argument establishing a multiplicative Chow--K\"unneth decomposition for cubic hypersurfaces \cite{Diaz}, \cite{FLV2}.

Using Theorem \ref{main}, we can prove a result concerning the {\em tautological ring\/}, which is a certain subring of the Chow ring of powers of $Y$:

\begin{nonumberingc}[=Corollary \ref{cor1}] Let $Y\subset\PP^{2g+1}$ be a smooth dimensionally transverse intersection of 2 quadrics. Let $m\in\NN$. Let
  \[ R^\ast(Y^m):=\Bigl\langle (p_i)^\ast(h), (p_{ij})^\ast(\Delta_Y)\Bigr\rangle\ \subset\ \ \ A^\ast(Y^m)   \]
  be the $\QQ$-subalgebra generated by pullbacks of the polarization $h\in A^1(Y)$ and pullbacks of the diagonal $\Delta_Y\in A^{2g-1}(Y\times Y)$. 
  The cycle class map induces injections
   \[ R^\ast(Y^m)\ \hookrightarrow\ H^\ast(Y^m,\QQ)\ \ \ \hbox{for\ all\ }m\in\NN\ .\]
   \end{nonumberingc}

 Corollary \ref{cor1} is inspired by results of Tavakol concerning the tautological ring of hyperelliptic curves \cite{Ta2}, \cite{Ta}. Loosely speaking, Corollary \ref{cor1} expresses that the Fano variety $Y$ behaves like a hyperelliptic curve from the point of view of intersection theory. This is in keeping with the known facts that $Y$
 is related to a genus $g$ hyperelliptic curve on the level of motives (Theorem \ref{YC} below), as well as on the level of derived categories (Remark \ref{DY} below).  
 
 It would be interesting to try and extend the results of the present paper. Do the same results hold for {\em all\/} complete intersections in projective space that have Hodge level 1 and dimension $>1$ (this is a very short list, cf. \cite[Table 1]{Rap}) ? Do the same results hold for intersections of 3 or more quadrics ?

   
%
%
   
 \vskip0.5cm

\begin{convention} In this article, the word {\sl variety\/} will refer to a reduced irreducible scheme of finite type over $\C$. A {\sl subvariety\/} is a (possibly reducible) reduced subscheme which is equidimensional. 

{\bf All Chow groups are with rational coefficients}: we will denote by $A_j(Y)$ the Chow group of $j$-dimensional cycles on $Y$ with $\QQ$-coefficients; for $Y$ smooth of dimension $n$ the notations $A_j(Y)$ and $A^{n-j}(Y)$ are used interchangeably. 
The notations $A^j_{hom}(Y)$ and $A^j_{AJ}(Y)$ will be used to indicate the subgroup of homologically trivial (resp. Abel--Jacobi trivial) cycles.
For a morphism $f\colon X\to Y$, we will write $\Gamma_f\in A_\ast(X\times Y)$ for the graph of $f$.

The contravariant category of Chow motives (i.e., pure motives with respect to rational equivalence as in \cite{Sc}, \cite{MNP}) will be denoted 
$\MM_{\rm rat}$.
\end{convention}

\section{Preliminaries}

\subsection{Intersections of 2 quadrics} 

Let $Y\subset\PP^{2g+1}$ be a smooth dimensionally transverse intersection of 2 quadrics.
The Hodge diamond of $Y$ looks like that of a curve (i.e. $h^{p,q}(Y)=0$ for $\vert p-q\vert>1$). This can be explained by the presence of a hyperelliptic curve
naturally associated to $Y$:

\begin{theorem}\label{YC} Let $Y\subset\PP^{2g+1}$ be a smooth dimensionally transverse intersection of 2 quadrics. There exists a hyperelliptic curve $C$ of genus $g$, such that
 there is an isomorphism of Chow motives
    \[ h(Y)\cong h^1(C)(1-g) \oplus \bigoplus_{j=0}^{2g-1}\one (-j)\ \ \ \hbox{in}\ \MM_{\rm rat}\ .\]
%
   \end{theorem}
   
  \begin{proof} That one can associate a hyperelliptic curve to $Y$ was already known to (and exploited by) Weil \cite{Weil}. The idea is as follows: 
  the two quadrics defining $Y$ span a pencil $\Gamma\cong\PP^1$. The choice of one of the two families of $g$-dimensional linear subspaces contained in a member of $\Gamma$ defines a double cover $C\to\Gamma$, which is ramified over the $2g+2$ points of $\Gamma$ corresponding to singular quadrics. Thus $C$ is a hyperelliptic curve of genus $g$. 
  Conversely, the Weierstrass points of any hyperelliptic curve determine a pencil of quadrics in $\PP^{2g+1}$, and the base locus of this pencil is a $Y$ as in the theorem.
  
 To prove the isomorphism of Chow motives, one first constructs a Chow--K\"unneth decomposition for $Y$ with the property that
   \begin{equation}\label{dck} h(Y)= h^{2g-1}(Y)\oplus   \bigoplus_{j=0}^{2g-1}\one (-j)\ \ \ \hbox{in}\ \MM_{\rm rat}\ \end{equation}
   (this is the Chow--K\"unneth decomposition which we will show to be multiplicative, cf. the proof of Theorem \ref{main} below).
   This is a standard construction: Letting $h\in A^1(Y)$ denote a hyperplane section (with respect to the embedding $Y\subset\PP^{2g+1}$), we consider
    \begin{equation}\label{ck}  \begin{split}  
                              \pi^{2j}_Y&:= {1\over 4}\, h^{2g-1-j}\times h^j\ \ \ \ \ (j=0,\ldots, 2g-1)\ ,\\
                                \pi^{2g-1}_Y&:= \Delta_Y-\sum_{j=0}^{2g-1} \pi^{2j}_Y\ \ \ \ \ \ \in\ A^{2g-1}(Y\times Y)\ .\\
                                \end{split}\end{equation}

  The existence of an isomorphism of homological motives
   \begin{equation}\label{isoh} h^{2g-1}(Y)\cong h^1(C)(1-g)\ \ \ \hbox{in}\ \MM_{\rm hom} \end{equation}
   is proven by Reid in his thesis \cite[Theorem 4.8]{Reid}. The idea is to use the variety $F=F_{g-1}(Y)$ parametrizing
 $g-1$-dimensional linear subspaces contained in $Y$. This is a $g$-dimensional abelian variety (cf. Theorem \ref{YF} below), and the isomorphism \eqref{isoh} is
 obtained by composing the isomorphisms
    \[ h^{2g-1}(Y)\ \xrightarrow{\cong} h^1(F)(1-g)  \ \xrightarrow{\cong}\ h^1(C)(1-g)\ ,\]
    where the first isomorphism is induced by the universal linear subspace (cf. Theorem \ref{YF} below), and the second isomorphism depends on the choice of a base point, cf. \cite[Proposition 4.2]{Reid}.

   As both sides of \eqref{isoh} are Kimura finite-dimensional (for $Y$ this is Lemma \ref{fd} below), this can be upgraded to an isomorphism of Chow motives
   \[ h^{2g-1}(Y)\cong h^1(C)(1-g)\ \ \ \hbox{in}\ \MM_{\rm rat}\ .\]  
   Combined with the decomposition \eqref{dck}, this proves the theorem.
   
  \begin{lemma}\label{fd} Let $Y\subset\PP^{2g+1}$ be a smooth dimensionally transverse intersection of 2 quadrics. Then $Y$ has finite-dimensional motive, in the sense of Kimura \cite{Kim}.
  \end{lemma}
  
  \begin{proof} It is known that
    \[ A^j(Y)=\QQ\ \ \ \forall\ j> g \]
    \cite{Ot} (cf. also \cite[Theorem 1.2]{BT} for a derived category argument). Using the Bloch--Srinivas argument \cite{BS}, this implies that
    \[ A^j_{AJ}(Y)=0\ \ \ \forall j\ .\]
    The lemma now follows from \cite[Theorem 4]{43}.
    \end{proof}
    
   \end{proof}

\begin{remark} Intersections of 2 quadrics in $\PP^5$ are one of the 17 families of Fano threefolds of Picard number one \cite{IP}; they are exactly the Fano threefolds of Picard number one that are {\em del Pezzo\/} (i.e. the canonical divisor $K_Y$ is $-2 H$, where $H$ is a positive generator of $\pic(Y)$) and of degree $4$ (i.e. $H^3=4$).
\end{remark}

\begin{remark}\label{DY} To give some context, we remark that the relation between $Y$ and $C$ of Theorem \ref{YC} also holds on the level of derived categories: there is a semi-orthogonal decomposition of the derived category
 \[ D^b(Y)=    \Bigl\langle \OO_Y(-2g+3),\ldots,\OO_Y(-1),\OO_Y,D^b(C)\Bigr\rangle\ \]
\cite{BO} (alternative proofs can be found in \cite{Ku0} and \cite{Add}). We will not use this fact.
 \end{remark}

%
%
%
%
%
 
 \subsection{Fano varieties of linear subspaces}
 
 \begin{theorem}\label{YF} Let $Y\subset\PP^{2g+1}$ be a smooth dimensionally transverse intersection of 2 quadrics. Let $F=F_{g-1}(Y)$ be the variety parametrizing
 $g-1$-dimensional linear subspaces contained in $Y$. Then $F$ is a smooth projective variety of dimension $g$, isomorphic to the Jacobian of the curve $C$ of Theorem \ref{YC}. The incidence correspondence $P\subset F\times Y$ induces an isomorphism of Chow motives
   \[  {}^t P\colon\ \ h^{2g-1}(Y)\ \xrightarrow{\cong}\ h^1(F)(1-g)\ \ \ \hbox{in}\ \MM_{\rm rat}\ .\]
    \end{theorem}
    
    \begin{proof} That $F$ is isomorphic to the Jacobian of $C$ is proven by Reid \cite[Theorem 4.8 and Theorem 4.14]{Reid} (cf. also \cite{DR} and \cite{Tyu} and \cite{Dona}).
  The isomorphism on the level of homological motives is \cite[Theorem 4.14]{Reid}; one then upgrades to Chow motives as above using Kimura finite-dimensionality.
%
      \end{proof}

 \subsection{MCK decomposition}
\label{ss:mck}

\begin{definition}[Murre \cite{Mur}] Let $X$ be a smooth projective variety of dimension $n$. We say that $X$ has a {\em CK decomposition\/} if there exists a decomposition of the diagonal
   \[ \Delta_X= \pi^0_X+ \pi^1_X+\cdots +\pi_X^{2n}\ \ \ \hbox{in}\ A^n(X\times X)\ ,\]
  such that the $\pi^i_X$ are mutually orthogonal idempotents and $(\pi_X^i)_\ast H^\ast(X,\QQ)= H^i(X,\QQ)$.
  
  (NB: ``CK decomposition'' is shorthand for ``Chow--K\"unneth decomposition''.)
\end{definition}

\begin{remark} The existence of a CK decomposition for any smooth projective variety is part of Murre's conjectures \cite{Mur}, \cite{J4}. 
\end{remark}

\begin{definition}[Shen--Vial \cite{SV}] Let $X$ be a smooth projective variety of dimension $n$. Let $\Delta_X^{sm}\in A^{2n}(X\times X\times X)$ be the class of the small diagonal
  \[ \Delta_X^{sm}:=\bigl\{ (x,x,x)\ \vert\ x\in X\bigr\}\ \subset\ X\times X\times X\ .\]
  An {\em MCK decomposition\/} is a CK decomposition $\{\pi_X^i\}$ of $X$ that is {\em multiplicative\/}, i.e. it satisfies
  \begin{equation}\label{vani} \pi_X^k\circ \Delta_X^{sm}\circ (\pi_X^i\times \pi_X^j)=0\ \ \ \hbox{in}\ A^{2n}(X\times X\times X)\ \ \ \hbox{for\ all\ }i+j\not=k\ .\end{equation}
  Here $\pi^i_X\times \pi^j_X$ is by definition $(p_{13})^\ast(\pi^i_X)\cdot (p_{24})^\ast(\pi^j_X)\in A^{2n}(X^4)$, where $p_{rs}\colon X^4\to X^2$ denotes projection on $r$th and $s$th factors.
  
 (NB: ``MCK decomposition'' is shorthand for ``multiplicative Chow--K\"unneth decomposition''.) 
  
  \end{definition}
  
  \begin{remark}\label{rem:mck} Note that the vanishing \eqref{vani} is always true modulo homological equivalence; this is because the cup product in cohomology respects the grading.
  
    The small diagonal (seen as a correspondence from $X\times X$ to $X$) induces the {\em multiplication morphism\/}
    \[ \Delta_X^{sm}\colon\ \  h(X)\otimes h(X)\ \to\ h(X)\ \ \ \hbox{in}\ \MM_{\rm rat}\ .\]
 Suppose $X$ has a CK decomposition
  \[ h(X)=\bigoplus_{i=0}^{2n} h^i(X)\ \ \ \hbox{in}\ \MM_{\rm rat}\ .\]
  By definition, this decomposition is multiplicative if for any $i,j$ the composition
  \[ h^i(X)\otimes h^j(X)\ \to\ h(X)\otimes h(X)\ \xrightarrow{\Delta_X^{sm}}\ h(X)\ \ \ \hbox{in}\ \MM_{\rm rat}\]
  factors through $h^{i+j}(X)$.
  
  If $X$ has an MCK decomposition, then setting
    \[ A^i_{(j)}(X):= (\pi_X^{2i-j})_\ast A^i(X) \ ,\]
    one obtains a bigraded ring structure on the Chow ring: that is, the intersection product sends $A^i_{(j)}(X)\otimes A^{i^\prime}_{(j^\prime)}(X) $ to  $A^{i+i^\prime}_{(j+j^\prime)}(X)$.
    
      It is expected that for any $X$ with an MCK decomposition, one has
    \[ A^i_{(j)}(X)\stackrel{??}{=}0\ \ \ \hbox{for}\ j<0\ ,\ \ \ A^i_{(0)}(X)\cap A^i_{hom}(X)\stackrel{??}{=}0\ ;\]
    this is related to Murre's conjectures B and D, that have been formulated for any CK decomposition \cite{Mur}.

  The property of having an MCK decomposition is severely restrictive, and is closely related to Beauville's ``splitting property' conjecture'' \cite{Beau3}. 
  To give an idea: hyperelliptic curves have an MCK decomposition \cite[Example 8.16]{SV}, but the very general curve of genus $\ge 3$ does not have an MCK decomposition \cite[Example 2.3]{FLV2}. As for surfaces: a smooth quartic in $\PP^3$ has an MCK decomposition, but a very general surface of degree $ \ge 7$ in $\PP^3$ should not have an MCK decomposition \cite[Proposition 3.4]{FLV2}.
For more detailed discussion, and examples of varieties with an MCK decomposition, we refer to \cite[Section 8]{SV}, as well as \cite{V6}, \cite{SV2}, \cite{FTV}, \cite{37}, \cite{38}, \cite{39}, \cite{40}, \cite{44}, \cite{46}, \cite{FLV2}, \cite{g8}.
   \end{remark}

 \section{The Franchetta property}
 
 \subsection{Definition}
 
 \begin{definition} Let $\XX\to B$ be a smooth projective morphism, where $\XX, B$ are smooth quasi-projective varieties. We say that $\XX\to B$ has the {\em Franchetta property in codimension $j$\/} if the following holds: for every $\Gamma\in A^j(\XX)$ such that the restriction $\Gamma\vert_{X_b}$ is homologically trivial for the very general $b\in B$, the restriction $\Gamma\vert_b$ is zero in $A^j(X_b)$ for all $b\in B$.
 
 We say that $\XX\to B$ has the {\em Franchetta property\/} if $\XX\to B$ has the Franchetta property in codimension $j$ for all $j$.
 \end{definition}
 
 This property is studied in \cite{PSY}, \cite{BL}, \cite{FLV}, \cite{FLV3}.
 
 \begin{definition} Given a family $\XX\to B$ as above, with $X:=X_b$ a fiber, we write
   \[ GDA^j_B(X):=\ima\Bigl( A^j(\XX)\to A^j(X)\Bigr) \]
   for the subgroup of {\em generically defined cycles}. 
   (In a context where it is clear to which family we are referring, the index $B$ will sometimes be suppressed from the notation.)
  \end{definition}
  
  With this notation, the Franchetta property amounts to saying that $GDA^\ast(X)$ injects into cohomology, under the cycle class map. 
 
\subsection{The families}

  \begin{notation}\label{not} Let 
   \[  \bar{B}:=\PP H^0(\PP^{2g+1},\OO_{\PP^{2g+1}}(2)^{\oplus 2}) \cong \PP^r\ ,\] 
   and let $B\subset\bar{B}$
   denote the Zariski open parametrizing smooth dimensionally transverse complete intersections. Let 
   \[ \YY\ \to\ B \]
   denote the universal family of smooth $2g-1$-dimensional complete intersections. 
  \end{notation}
  
  \begin{notation}\label{not2} Let $B$ be as in Notation \ref{not}, and let $\FF\to B$ denote the universal family of Fano varieties of $g-1$-dimensional linear subspaces contained in the intersection of 2 quadrics. Here $\FF\subset \Gr(g,2g+2)\times B$ is defined as
    \[ \FF:= \Bigl\{  (L, b)\ \big\vert L\subset Y_b \Bigr\}\ \ \ \subset \ \Gr(g,2g+2)\times B  \ .\]
       \end{notation}
  
  \begin{lemma}\label{F} The variety $\FF$ is smooth.
   \end{lemma}
   
   \begin{proof} One way to see this is by observing that $\FF\to B$ is the relative Jacobian of a smooth relative curve (Theorem \ref{YF}). One can also see the smoothness more directly as follows:
   
   Let $\bar{\FF}\supset \FF$ denote the Zariski closure of $\FF$ in $\Gr(g,2g+2)\times \bar{B}$.
   Let $S\to \Gr(g,2g+2)$ be the tautological rank $g$ vector bundle. An element $b\in B$ determines a section of the vector bundle $E:=\sym^2 S\oplus \sym^2 S$, and the zero locus of this section is the Fano variety $F(Y_b)$.
   
   A point $L\in\Gr(g,2g+2)$ imposes $2 {{g+1}\choose{2}}=\rank E$ conditions on sections of $E$, and so the projection
     \[  \bar{\FF}\ \to\ \Gr(g,2g+2) \]
     is a $\PP^{s}$-bundle with $s=r-2  {{g+1}\choose{2}}$. It follows that $\bar{\FF}$ is smooth. As $\FF\subset\bar{\FF}$ is a Zariski open, $\FF$ is also smooth.       
   \end{proof}

  \subsection{Franchetta for $Y$ and $Y\times Y$}
 
 \begin{proposition}\label{Fr} Let $\YY\to B$ and $\FF\to B$ be the universal families as in Notations \ref{not} and \ref{not2}.
 The following families have the Franchetta property:
 
 \noindent
 (\rom1)
 The family $\YY\to B$; 
 
 \noindent
 (\rom2) the family
 $\YY\times_B \YY\to B$.
%
  \end{proposition}
 
\begin{proof} 

\noindent
(\rom1)
Let $\bar{\YY}\to\bar{B}$ denote the Zariski closure. This is the universal family of all (possibly singular and degenerate) intersections of 2 quadrics in $\PP^{2g+1}$.
Since every point in $\PP^{2g+1}$ imposes one condition on each quadric, the projection $\pi\colon\bar{\YY}\to \PP^{2g+1}$ is a $\PP^{r-2}$-bundle.  
 
 Now given $\Gamma\in A^j({\YY})$, let $\bar{\Gamma}\in A^{j}(\bar{\YY})$ be a cycle restricting to $\Gamma$. The projective bundle
		formula yields
		\begin{equation}\label{pb}  \bar{\Gamma}=  \pi^\ast(a_{j}) + \pi^\ast(a_{j-1})\cdot \xi   +       \cdots + \pi^\ast(a_0)\cdot \xi^j    \ \ \  \hbox{in}\
		A^{j}(\bar{\YY}) ,\end{equation}
		where $\xi\in A^1(\bar{\YY})$ is relatively ample with respect to $\pi$ and
		$a_j\in A^{j}(\PP^{2g+1})$.
		Let $h\in A^1(\bar{B})$ be a hyperplane section and let $\nu\colon \bar{\YY}
		\to \bar{B}$ denote the projection. We have
		\[  \nu^\ast(h)= c\, \xi+\pi^\ast(z)\ \ \ \hbox{in}\ A^{1}(\bar{\YY})\ ,\] 
		for some $c \in\QQ$ and $z\in A^{1}(\PP^{2g+1})$. It
		is readily checked that the constant $c$ is non-zero. (Indeed, let us assume for a moment $c$
		were zero. Then we would have $\nu^\ast(h^{r})=\pi^\ast(z^{r})$
		in $A^{r}(\bar{\YY})$. But the right-hand side is zero, since
		$r:=\dim\bar{B}>\dim \PP^{2g+1}=2g+1$, while the left-hand side is non-zero.)	
		
		The constant $c$ being non-zero, we can write
		\[ \xi= \pi^\ast(z)+\nu^\ast(d)\ \ \  \hbox{in}\  A^{1}(\bar{\YY}) ,\]
		where $z\in A^{1}(\PP^{2g+1})$ and $d\in A^1(\bar{B})$ are non-zero elements.
		The restriction of $\nu^\ast(d)$ to a fiber $Y_b$ is zero,
		and so equality \eqref{pb} implies that
		\[  {\Gamma}\vert_{Y_b}=  a_j^\prime\vert_{Y_b}\ \ \  \hbox{in}\
		A^j(Y_b) ,\]
		for some $a_j^\prime\in A^{j}(\PP^{2g+1})$. Since $A^{j}(\PP^{2g+1})\cong\QQ$, this proves the Franchetta property for $\YY\to B$.
		
\medskip
\noindent
(\rom2) Let us now consider the family $\YY\times_B \YY\to B$. As above, we have the Zariski closure $B\subset \bar{B}\cong\PP^r$, and the induced family $\bar{\YY}\times_{\bar{B}}\bar{\YY}\to\bar{B}$. Since 2 distinct points in $\PP^{2g+1}$ impose 2 independent conditions on each quadric, the morphism
  \[ (\pi,\pi)\colon\ \ \bar{\YY}\times_{\bar{B}}\bar{\YY}\ \to\ \PP^{2g+1}\times \PP^{2g+1} \]
  has the structure of a $\PP^{r- 4}$-bundle over $(\PP^{2g+1}\times \PP^{2g+1} )\setminus \Delta_{\PP^{2g+1}}$, and a $\PP^{r-2}$-bundle over $  \Delta_{\PP^{2g+1}}$.
  That is, $(\pi,\pi)$ is a {\em stratified projective bundle\/}, in the sense of \cite[Definition 5.1]{FLV}. Applying \cite[Proposition 5.2]{FLV},
one thus find that there is equality
    \begin{equation}\label{square} \begin{split} \ima\Bigl(  A^\ast(\bar{\YY}\times_{\bar{B}}\bar{\YY})\to A^\ast(Y_b\times Y_b)\Bigr)      = \ima\Bigl(  A^\ast(\PP^{2g+1}\times \PP^{2g+1})\to  A^\ast(Y_b\times Y_b)\Bigr)&\\  
    + \Delta_\ast \ima\Bigl( A^\ast(\PP^{2g+1})\to &A^\ast(Y_b)\Bigr)\ ,\\
    \end{split}\end{equation}
 where $\Delta\colon Y_b\to Y_b\times Y_b$ denotes the diagonal embedding.
 
 Clearly one has equality
  \[\ima\Bigl( A^\ast(\PP^{2g+1})\to A^\ast(Y_b)\Bigr) = \langle h \rangle\ .\] 
  Moreover, it is well-known that
   \[ A^\ast(\PP^{2g+1}\times \PP^{2g+1})=A^\ast(\PP^{2g+1})\otimes A^\ast(\PP^{2g+1})\ .\]
   In view of these two observations, equality \eqref{square} reduces to
   \[   \ima\Bigl(  A^\ast(\bar{\YY}\times_{\bar{B}}\bar{\YY})\to A^\ast(Y_b\times Y_b)\Bigr)      =  \bigl\langle h, \Delta_{Y_b} \bigr\rangle\ .\]
   It now remains to prove that for any smooth $Y_b$ the $\QQ$-subalgebra 
     \[  R^\ast(Y_b\times Y_b):=\bigl\langle h, \Delta_{Y_b} \bigr\rangle \ \  \subset\  A^\ast(Y_b\times Y_b)\]
   injects into cohomology via the cycle class map. 
   
   To this end, let $\Lambda$ denote the moduli stack of smooth dimensionally transverse intersections of 2 quadrics in $\PP^{g+1}$, and let $\YY\to \Lambda$ denote the universal family (this exists as an algebraic stack). The $\QQ$-subalgebra $R^\ast=\langle \Delta_{Y_b}, K_{Y_b}\rangle$ is contained in (and hence equal to) $GDA^\ast_\Lambda(Y_b\times Y_b)$, and so we are reduced to proving the Franchetta property for $\YY\times_\Lambda \YY$. 
   
Let $\HH_g\subset\MM_g$ denote the image of the period map $\Lambda\to \MM_g$ (this $\HH_g$ is the locus of jacobians of hyperelliptic curves of genus $g$).
Because any hyperelliptic curve of genus $g$ is attained by the construction of Theorem \ref{YC}, and conversely $Y$ is determined by the hyperelliptic curve $C$, the morphism from the moduli stack $\Lambda$ to $\HH_g$ is injective, cf. \cite[Section 3.6]{Deb}.
The injectivity means that we may as well consider $Y$ as a fiber of the family $\YY\to \HH_g$, and prove the Franchetta property for $\YY\times_{\HH_g} \YY$.  
 The isomorphism of motives
  \[ h(Y)\cong h^1(C)(1-g)\oplus \bigoplus \one(\ast)\ \ \ \hbox{in}\ \MM_{\rm rat} \]
  induces a split injection of Chow groups
   \[ A^j (Y\times Y)\ \hookrightarrow\ A^{j+2-2g}(C\times C ) \oplus \bigoplus A^\ast(C)\ ,\]
   and this injection preserves generically defined cycles 
   with respect to $\HH_g$ (indeed, both $P$ and the isomorphism $h^1(F)\cong h^1(C)$ are generically defined). This means that there is an injection
   \[  GDA^j_{\HH_g}(Y\times Y)\ \hookrightarrow\ GDA^{j+2-2g}_{\HH_g}(C\times C )\oplus \bigoplus GDA^\ast_{\HH_g}(C)\ .\]
   The Franchetta property for $\YY\times_\Lambda \YY$
   thus reduces to the following lemma:
   
    \begin{lemma}\label{Hg} Let  $\CC\to\HH_g$ denote the universal hyperelliptic curve of genus $g$. The family $\CC\times_{\HH_g}\CC\to \HH_g$
   has the Franchetta property.
   \end{lemma}   
   
   To prove the lemma, we use that the general hyperelliptic curve $C$ of genus $g$ can be realized by a degree $2g+2$ polynomial $x_2^2=f(x_0,x_1)$ in the weighted projective plane
   $\PP:=\PP(1,1,g+1)$. This $\PP$ can be seen as the cone over the rational normal curve of degree $g+1$ in $\PP^{g+1}$, and $C$ corresponds to a quadric section of this cone avoiding the singular point $[0:0:1]\in\PP$. Let $B$ denote the parameter space parametrizing smooth quadric sections of $\PP$, so $B$ is an open in some $\bar{B}\cong\PP^s$, and let $\CC\to B$ denote the universal family.
   We have an injection
     \[ GDA^\ast_{\HH_g}(C\times C)\ \hookrightarrow\ GDA^\ast_{B}(C\times C) \ ,\]
   and so it suffices to prove the Franchetta property for $\CC\times_B \CC\to B$. This is easy: the projection
     \[ \bar{\CC}\times_{\bar{B}}\bar{\CC}\ \to\ \PP\times\PP \]
     is a $\PP^{s-2}$-fibration over $(\PP\times\PP)\setminus \Delta_\PP$ and a $\PP^{s-1}$-fibration over $\Delta_\PP$. The stratified projective bundle argument \cite[Proposition 5.2]{FLV} then yields the equality
     \[ GDA^\ast_{B}(C\times C)= \langle h_C,\Delta_C\rangle\ \]
     (where $h_C\in A^1(C)$ is the restriction of a hyperplane in $\PP^{g+1}$).
     That the right-hand side injects into cohomology for any hyperelliptic curve $C$ of genus $g\ge 2$ was proven by Tavakol \cite{Ta}; alternatively one can see this directly in the following way. For codimension 1, the diagonal in cohomology is not a combination of decomposable cycles (indeed, $\Delta_C$ acts as the identity on $H^1(C,\QQ)$), and so injectivity reduces to the Franchetta property for $\CC\to B$. In codimension 2, there is equality
     \[ \Delta_C\cdot (p_1)^\ast(h) = \lambda\, {}^t \Gamma_\tau\circ \Gamma_\tau\ \ \ \hbox{in}\ A^{2}(C\times C)\ ,\]
     where $\lambda\in\QQ$ and $\tau\colon C\hookrightarrow \PP$ denotes the inclusion morphism (NB: $\PP$ is a projective quotient variety and so the formalism of correspondences with rational coefficients still makes sense for $\PP$, as noted in \cite[Example 16.1.13]{F}). Using Lieberman's lemma, we then find that
     \[ \Delta_C\cdot (p_1)^\ast(h) =    \lambda\, (\tau\times\tau)^\ast (\Delta_\PP)= \lambda^\prime\, \sum_{j=0}^2 h_C^j\times h_C^{2-j}= \lambda^\prime \, h_C\times h_C\   \ \ 
     \hbox{in}\ A^{2}(C\times C)\ ,\]     
     and so the injectivity of $ \langle h_C,\Delta_C\rangle\cap A^2(C\times C)$ into $H^4(C\times C,\QQ)$ again follows.     
     
     The lemma, and hence the proposition, are now proven.
\end{proof}

\begin{corollary}\label{YY} Let $Y\subset\PP^{2g+1}$ be a smooth dimensionally transverse intersection of 2 quadrics. 
Then
  \[  GDA_{}^\ast(Y\times Y) = \Bigl\langle (p_1)^\ast GDA_{}^\ast(Y),  (p_2)^\ast GDA_{}^\ast(Y)\Bigr\rangle \oplus \QQ[\Delta_Y]\ .\]
  In particular, there is equality
     \begin{equation}\label{hyp} \Delta_{Y}\cdot (p_j)^\ast(h) =\sum_{i=1}^{2g}  a_i\,  h^i\times h^{2g-i}\ \ \ \hbox{in}\ A^{2g}(Y\times Y) \ \ \ (j=1,2)\ ,\end{equation}
   for some $a_i\in\QQ$ (actually, all $a_i$ are $1\over 4$). 
      \end{corollary}
      
      \begin{proof} Equality \eqref{hyp} follows immediately from the first equality, since 
        \[ \Delta_{Y}\cdot (p_j)^\ast(h)\in GDA^{2g}(Y\times Y)\ .\]
        (To find the values of $a_i$, one observes that $(\Delta_Y - {1\over 4}\sum_{i=0}^{2g-1} h^i\times h^{2g-1-i})\cdot (p_j)^\ast(h)$ acts as zero on cohomology.)
      
      To prove the first equality, we recall from the proof of Proposition \ref{Fr}(\rom2) the equality
      \[ GDA^\ast(Y\times Y)=\langle K_Y,\Delta_Y\rangle = GDA^\ast_{\Lambda}(Y\times Y)\ ,\]
      and the split injection
      \begin{equation}\label{sin} GDA^j_\Lambda(Y\times Y)\ \hookrightarrow\ GDA^{j+2-2g}_{\HH_g}(C\times C)\oplus \bigoplus GDA^\ast(C)\ .\end{equation}
      We have seen that for any $j+2-2g\not=1$, the group $GDA^{j+2-2g}_{\HH_g}(C\times C)$ is {\em decomposable\/}, i.e.
      \[        GDA^{j+2-2g}_{\HH_g}(C\times C) = \Bigl\langle (p_1)^\ast GDA_{\HH_g}(C),     (p_2)^\ast GDA_{\HH_g}(C)\Bigr\rangle\ \ \ \forall j+2-2g\not=1\ .\]
      The left-inverse of the split injection \eqref{sin} sends decomposable cycles to decomposable cycles, and so we find that
      \[ GDA^j_{\Lambda}(Y\times Y) = \Bigl\langle (p_1)^\ast GDA_{\Lambda}(Y),     (p_2)^\ast GDA_{\Lambda}(Y)\Bigr\rangle\ \ \ \forall j\not=2g-1\ .\]     
%
Because we know that $\Delta_Y$ is linearly independent from the decomposable cycles (indeed, $\Delta_Y$ acts as the identity on $H^{2g-1}(Y,\QQ)$), this proves the corollary.       
      \end{proof}

\begin{remark}\label{fp} 
The equality \eqref{hyp} is remarkable, for the following reason. In case $Y\subset\PP^{m}$ is a smooth {\em hypersurface\/} (of any degree), equality \eqref{hyp}
is true for $Y$, as follows from the excess intersection formula (or from the argument of Lemma \ref{Hg}). On the other hand, in case $Y\subset\PP^m$ is a complete intersection of codimension at least 2,
in general there is {\em no equality\/} of the form \eqref{hyp}. Indeed, let $C$ be a very general curve of genus $g\ge 4$. The Faber--Pandharipande cycle
  \[  FP(C):= \Delta_C\cdot (p_j)^\ast(K_C) - {1\over 2g-2} K_C\times K_C\ \ \ \in A^2(C\times C)\ \ \ \ \ (j=1,2) \]
  is homologically trivial but non-zero in $A^2(C\times C)$ \cite{GG}, \cite{Yin0} (this cycle $FP(C)$ is the ``interesting 0-cycle'' in the title of \cite{GG}). In particular, for the very general complete intersection $Y\subset\PP^3$ of bidegree $(2,3)$, the cycle
  \[ FP(C):=  \Delta_C\cdot (p_j)^\ast(h) - {1\over 2g-2} h \times h\ \ \ \in A^2(C\times C) \]
  is homologically trivial but non-zero, and so there cannot exist an equality of the form \eqref{hyp} for $Y$.
%
\end{remark}

\subsection{Franchetta for $F\times Y$ (the case $g=2$)} We now restrict to Fano threefolds $Y\subset\PP^5$ that are intersections of 2 quadrics, and their Fano surfaces of lines 
$F=F(Y)$. In this setting, we have the following Franchetta property (the result is probably true for intersections of 2 quadrics $Y\subset\PP^{2g+1}$ for arbitrary $g$, but the argument is easier in case $g=2$, which is the only case we will need in this paper).

\begin{proposition}\label{Fr2} Assume $g=2$, and let $\YY\to B$ and $\FF\to B$ be the universal families as in Notations \ref{not} and \ref{not2} (i.e. $\YY\to B$ is a family of Fano threefolds and $\FF\to B$ is a family of Fano surfaces). The family
$\FF\times_B \YY\to B$ has the Franchetta property.
\end{proposition}

\begin{proof} As a warm-up, let us first prove the following:

\begin{proposition}\label{Fr3} The family $\FF\to B$ has the Franchetta property. More precisely,
  \[ GDA^\ast(F)=\ima\Bigl( A^\ast\bigl(\Gr(2,6)\bigr)\to A^\ast(F)\Bigr)\ \]
 injects into cohomology under the cycle class map.
 \end{proposition}
 
 \begin{proof}
As we have already seen in the proof of Lemma \ref{F}, 
$\bar{\FF}\to \Gr(2,6)$ is a projective bundle (this is just because a line imposes exactly 3 conditions on quadrics [indeed, a quadric containing 3 distinct points on the line has to contain the line], and so a line imposes 6 conditions on $\bar{B}$). 
Applying the projective bundle formula, and reasoning as in the proof of Proposition \ref{Fr} (or directly applying
\cite[Proposition 5.2]{FLV}), this implies the equality
   \[ GDA^\ast(F)=\ima\Bigl( A^\ast\bigl(\Gr(2,6)\bigr)\to A^\ast(F)\Bigr)\ .\]
  Since $F$ is an abelian surface, its motive decomposes
    \[ h(F)= \one \oplus h^1(F) \oplus \sym^2 h^1(F) \oplus h^1(F)(-1) \oplus \one(-2)\ \ \ \hbox{in}\ \MM_{\rm rat}\ .\]
    We know (Theorem \ref{YF}) that $h^1(F)$ is isomorphic to $h^3(Y)(1)$. This induces a split injection of motives
    \[ h(F)\ \hookrightarrow\ h(Y\times Y)(2) \oplus h(Y)(1) \oplus h(Y) \oplus \one \oplus \one(-2)\ ,\]
    and in particular a split injection of Chow groups
    \begin{equation}\label{split} A^j(F)\ \hookrightarrow\ A^{j+2}(Y\times Y)\oplus A^{j+1}(Y)\oplus A^j(Y)\oplus \QQ^2\ .\end{equation}
    Since the decomposition of $h(F)$ and the Abel--Jacobi isomorphism $h^1(F)\cong h^3(Y)(1)$ are generically defined, the injection \eqref{split}
    preserves generically defined cycles. The Franchetta property for $\FF\to B$ then follows from the Franchetta property for $\YY\times_B \YY\to B$ and for $\YY\to B$.
    \end{proof}
   
We now start the proof of Proposition \ref{Fr2}.
Let us write 
      \[ \Gamma\ \ \subset\ \Gr(2,6)\times\PP^5\] 
      for the universal line.
Given a line $\ell\subset\PP^5$ and a point $x\in\PP^5$, $x\not\in\ell$, there exist quadrics containing $\ell$ and avoiding $x$ (indeed, one may simply take degenerate quadrics consisting of a double $\PP^4$ containing $\ell$ and avoiding $x$). This means that the pair $(\ell,x)$ with  $x\not\in\ell$ imposes
exactly 8 conditions on $\bar{B}$.
 It follows that
the projection
    \[ \pi\colon\ \  \bar{\FF}\times_{\bar{B}} \bar{\YY}\ \to\  \Gr(2,6)\times \PP^5 \]
    has the structure of a $\PP^{r-8}$-bundle over $(\Gr(2,6)\times \PP^5)\setminus\Gamma$, and a $\PP^{r-6}$-bundle over $\Gamma$ 
    That is, $\pi$ is a {\em stratified projective bundle\/}, in the sense of \cite{FLV}. Applying \cite[Proposition 5.2]{FLV}, we find that 
    \[ \begin{split}GDA^\ast(F\times Y):=\ima\Bigl( A_\ast( \bar{\FF}\times_{\bar{B}} \bar{\YY})\to A^\ast(F\times Y)\Bigr) =     \ima\Bigl( A^\ast\bigl(\Gr(2,6)\times\PP^5\bigr)\to A^\ast(F\times Y)
        \Bigr)&\\ +   
              \tau_\ast \ima\Bigl( A^\ast(\Gamma)\to A^\ast(P)\Bigr)&\ ,\\
              \end{split}\]
              where $\tau\colon P\to F\times Y$ denotes the inclusion morphism. Since $\Gamma$ is a $\PP^1$-bundle over $\Gr(2,6)$ with $p^\ast(h)$ relatively ample (where $h\in A^1(\PP^5)$ is ample), pullback induces a surjection 
                \[   A^\ast\bigl(\Gr(2,6)\times\PP^5\bigr)\ \twoheadrightarrow \ A^\ast(\Gamma)\ ,\] 
                and the above simplifies to
          \[ GDA^\ast(F\times Y) =  \Bigl\langle    \ima\Bigl( A_\ast(\Gr(2,6)\times\PP^5)\to A^\ast(F\times Y)\Bigr)    , P\Bigr\rangle\ .\]        
       Since $\PP^5$ has (trivial Chow groups and hence) the Chow--K\"unneth property, this further reduces to
        \begin{equation}\label{gd} GDA^\ast(F\times Y) \ \subset\ \Bigl\langle  (p_F)^\ast GDA^\ast(F)  ,
                  (p_Y)^\ast GDA^\ast(Y), P \Bigr\rangle\ .\end{equation}
       
  Let us now have a closer look at the right-hand side of \eqref{gd}. We define the graded $\QQ$-vector space
    \[ \begin{split}   R^\ast(F\times Y):=  \Bigl\langle  (p_{F})^\ast GDA^\ast(F) , 
                  (p_Y)^\ast GDA^\ast(Y) \Bigr\rangle\   + \QQ[P] + \QQ[ P\cdot (p_{F})^\ast(h_{F})]&\\
                          \subset\ A^\ast(F&\times Y)\ ,\\
                          \end{split}
                 \]
                  where $h_{F}$ is a generator of $GDA^1(F)$ (note that $GDA^1(F)$ is 1-dimensional by Proposition \ref{Fr3}). Here $R^\ast(F\times Y)$ is naturally graded by codimension, which means that $P$ is placed in $R^2(F\times Y)$ and $P\cdot (p_{F})^\ast(h_{F})$ is placed in $R^3(F\times Y)$.
                  
    Clearly, $R^\ast(F\times Y)$ is contained in the right-hand side of \eqref{gd}, i.e. there is an inclusion of $\QQ$-vector spaces
    \begin{equation}\label{subset} R^\ast(F\times Y)\ \ \subset \   \Bigl\langle  (p_{F})^\ast GDA^\ast(F)  ,
                  (p_Y)^\ast GDA^\ast(Y), P \Bigr\rangle\ .      \end{equation}
    The goal of the next three lemmata will be to show that the inclusion \eqref{subset} is actually an equality (and hence, in particular, $R^\ast(F\times Y)\subset A^\ast(F\times Y)$ is actually a $\QQ$-subalgebra).              
                      
  \begin{lemma}\label{lemma1}
                 \[  P\cdot  (p_Y)^\ast GDA^\ast(Y)\ \ \in\ R^\ast(F\times Y)\ .\]
                               \end{lemma}
                  
     \begin{proof} Clearly
  \[ GDA^\ast(Y)= \ima\bigl( A^\ast(\PP^5)\to A^\ast(Y)\bigr)=\langle h\rangle \ \] 
  (where $h\in A^1(Y)$ denotes the polarization).     
  
As we have seen above (equality \eqref{hyp}) there exist $a_j\in\QQ$ such that
    \begin{equation}\label{req}
        \Delta_Y\cdot (p_2)^\ast(h) = \sum_j a_j\, h^j\times h^{4-j}\ \ \ \hbox{in}\ A^4(Y\times Y)\ . \end{equation}

  Let $p$ and $q$ denote the projections from $P$ to $F$ resp. $Y$.
  Using equality \eqref{req}, we find that
    \[ \begin{split} P\cdot (p_Y)^\ast(h) &=  (p\times\Delta_Y)_\ast (q\times\Delta_Y)^\ast    \Bigl(\Delta_Y\cdot (p_2)^\ast(h)\Bigr)\\
                         &=   (p\times\Delta_Y)_\ast (q\times\Delta_Y)^\ast  \Bigl(      \sum_j a_j\, h^j\times h^{4-j}\Bigr)\\
                         &= \sum_j a_j\, (p_F)^\ast P^\ast(h^j)\cdot (p_Y)^\ast (h^{4-j})\ \ \ \ \hbox{in}\ A^3(F\times Y)\ .\\
                         \end{split}\]
     Since $P$ and $h$ are generically defined, $P^\ast(h^j)\in GDA^\ast(F)$, and so we get
     \[  P\cdot (p_Y)^\ast(h)    \ \ \in \  \Bigl\langle  (p_F)^\ast GDA^\ast(F) , 
                  (p_Y)^\ast GDA^\ast(Y) \Bigr\rangle\ .\]
    It follows that likewise
     \begin{equation}\label{eqlem}     P\cdot (p_Y)^\ast(h^i)    \ \ \in \  \Bigl\langle  (p_F)^\ast GDA^\ast(F) , 
                  (p_Y)^\ast GDA^\ast(Y) \Bigr\rangle\ \ \ \subset\ R^\ast(F\times Y)\ ,\end{equation}
                  which proves the lemma.
        \end{proof}             
                  
   \begin{lemma}\label{lemma2}
    \[  P\cdot (p_F)^\ast GDA^\ast(F) \ \ \in\  R^\ast(F\times Y)    \ .\]      
                   \end{lemma}
                   
     \begin{proof} We have seen above (Proposition \ref{Fr3}) that
     \[ GDA^\ast(F)=\ima\Bigl( A^\ast\bigl(\Gr(2,6)\bigr)\to A^\ast(F)\Bigr)\ ,\]
     and so $GDA^1(F)$ is $1$-dimensional, generated by $h_F$. Since $\FF\to B$ has the Franchetta property (Proposition \ref{Fr3}),
     $GDA^2(F)$ is also $1$-dimensional, generated by $h_F^2$. It is readily checked (cf. Sublemma \ref{sub} below) that $GDA^2(F)$ is also generated by     
     $c_2(Q)\vert_F$, where $Q$ is the universal quotient bundle on $\Gr(2,6)$. To prove the lemma, we thus need to check that
       \begin{equation}\label{inth} P\cdot (p_F)^\ast ( c_2(Q)\vert_F)\ \ \in\  R^4(F\times Y)\ .\end{equation}
     
     The morphism $p\colon P\to F$ being a $\PP^1$-bundle (with $q^\ast(h)$ being relatively ample), we find that
       \[   p^\ast (c_2(Q)\vert_F) = -   q^\ast(h^2) - q^\ast(h) p^\ast (c_1(Q)\vert_F)  \ \ \ \hbox{in}\ A^2(P)\ .\]
       Pushing forward under the closed immersion $P\hookrightarrow F\times Y$, this implies that
       \[  P\cdot (p_F)^\ast ( c_2(Q)\vert_F) =  - P\cdot (p_Y)^\ast (h^2) - P\cdot (p_Y)^\ast(h)\cdot (p_F)^\ast(h_F)\ \ \ \hbox{in}\ A^4(F\times Y)\ .\]
      Using equation \eqref{eqlem}, we see that the right-hand side is contained in $R^4(F\times Y)$, and so \eqref{inth} is proven. 
      
         \begin{sublemma}\label{sub} Let $Q$ be the universal quotient bundle on $\Gr(2,6)$. The restriction $c_2(Q)\vert_F$ is non-zero in $H^4(F,\QQ)$.
       \end{sublemma}
       
        \begin{proof} This is an application of Schubert calculus, for which we adhere to the notation of \cite[Chapter 14]{F}. As $Y$ is intersection of 2 quadrics $Y=Q\cap Q^\prime$, the class of 
 $F$ in $A^\ast(\Gr(2,6))$ is the intersection $F_1(Q)\cdot F_1(Q^\prime)$. Using \cite[Example 14.7.14]{F}, we thus find 
   \[  F= 16 \{ 2,1\}\cdot \{2,1\}=16 (2,4)\cdot (2,4)\ \ \ \hbox{in}\ A^6(\Gr(2,6))\ .\]
   Using Pieri's rule \cite[Proposition 14.6.1]{F}, it follows that
   \[ \begin{split} \deg (c_2(Q)\vert_F) &=  16 (2,4)\cdot (2,4)\cdot (2,5) \\
                                   &= 16 (2,4)\cdot \Bigl( (0,4) + (1,3)\Bigr) \\
                                   &=16\ .\\
                                   \end{split}\]
     \end{proof} 
     
     The sublemma, and hence the lemma, are proven.
     \end{proof}              
                   
    \begin{lemma}\label{lemma3}
     \[ P^2\ \ \in\  R^4(F\times Y)\ .\]
     \end{lemma}               
                  
             \begin{proof} The self-intersection formula \cite[Corollary 6.3]{F} gives the equality
               \[   P^2 = \tau_\ast c_2(N)\ \ \ \hbox{in}\ A^4(F\times Y)\ ,\]
               where $\tau\colon P\to F\times Y$ is the inclusion morphism and
               $N$ is the vector bundle on $P$ fitting into an exact sequence
                \begin{equation}\label{se} 0 \ \to\ TP\ \to\   \tau^\ast(TF\oplus TY)\ \to\ N \ \to\ 0\ .\end{equation}
             Since $p\colon P\to F$ is a $\PP^1$-bundle, the Chern classes of $P$ can be computed from the exact sequence
             \[ 0 \to\  T_p\ \to\ T_P\ \to\ p^\ast T_F\ \to\ 0\ ,\]
             where $T_p$ fits into
             the Euler sequence
             \[ 0 \ \to\   \OO_P\ \to\ p^\ast\EE(1)\ \to\ T_p\ \ \to\ 0\ .\]
             (Here $\EE$ is the restriction to $F$ of the tautological subbundle on $\Gr(2,6)$.)
             These exact sequences imply that the Chern classes $c_j(P)$ are in $\langle p^\ast c_j(Q)\vert_F, q^\ast(h)\rangle$. The exact sequence \eqref{se}
             implies that the same holds for the Chern classes $c_j(N)$. Pushing forward under $\tau$, this gives
             \[ P^2\ \ \in\ P\cdot \Bigl\langle  (p_F)^\ast GDA^\ast(F), (p_Y)^\ast GDA^\ast(Y)  \Bigr\rangle\ .\]
             Using Lemmata \ref{lemma1} and \ref{lemma2}, it follows that $P^2\in R^4(F\times Y)$.
                  \end{proof}

  Applying Lemmata \ref{lemma1}, \ref{lemma2} and \ref{lemma3}, we  conclude that (as promised) the inclusion \eqref{subset} is an equality, i.e.
            \[         \Bigl\langle  (p_{F})^\ast GDA^\ast(F)  ,
                  (p_Y)^\ast GDA^\ast(Y), P \Bigr\rangle =   R^\ast(F\times Y)\ .       \]
   Combining with equality \eqref{gd}, we obtain an equality
   \begin{equation}\label{GDR}  GDA^\ast(F\times Y) = R^\ast(F\times Y) \ .\end{equation}
   
   Using equality \eqref{GDR}, we are now in position to check that the Franchetta property holds for $\FF\times_B \YY\to B$.
   Recall that by definition
    \[  R^\ast(F\times Y):=   \Bigl\langle  (p_{F})^\ast GDA^\ast(F) , 
                  (p_Y)^\ast GDA^\ast(Y) \Bigr\rangle\   + \QQ[P] + \QQ[ P\cdot (p_{F})^\ast(h_{F})]\ .\]       
  As such, equality \eqref{GDR} implies that $GDA^j(F\times Y)$ is {\em decomposable\/} for $j\not=2,3$, i.e.
    \[  GDA^j(F\times Y) =    \Bigl\langle  (p_{F})^\ast GDA^\ast(F) , 
                  (p_Y)^\ast GDA^\ast(Y) \Bigr\rangle\ \ \ \forall j\not=2,3\ .\]
  The Franchetta property for $\FF\times_B \YY$ in codimension $j\not=2,3$ thus follows from the Franchetta property for $\FF\to B$ and for $\YY\to B$.
  
  For codimension $j=2$, we remark that the class of $P$ in $H^4(F\times Y,\QQ)$ cannot be decomposable: indeed, if it were, the transpose ${}^t P$ would act as zero on $H^3(Y,\QQ)$, contradicting the Abel--Jacobi isomorphism 
    \[  P^\ast\colon\ \  H^3(Y,\QQ)\ \xrightarrow{\cong}\ H^1(F,\QQ)\] 
    (Theorem \ref{YF}). 
Hence, $P$ and the decomposable part of $R^2(F\times Y)$ are linearly independent in cohomology, and the Franchetta property in codimension 2 again follows from that for the families $\FF\to B$ and for $\YY\to B$.
  
  Likewise, in codimension $j=3$, the class $P\cdot    (p_{F})^\ast(h_{F})$ is not decomposable (because combining the Abel--Jacobi isomorphism with hard Lefschetz, one obtains an isomorphism $H^3(Y,\QQ)\cong H^3(F,\QQ)$ which is induced by $P\cdot    (p_{F})^\ast(h_{F})$), and the Franchetta property in codimension 3 follows. The proposition is proven.                         
          \end{proof}

  \section{Main result}
  
  \begin{theorem}\label{main} Let $Y\subset\PP^{2g+1}$ be a smooth dimensionally transverse intersection of 2 quadrics. Then $Y$ has a multiplicative Chow--K\"unneth decomposition.
  \end{theorem}

 \begin{proof} In case $g=1$, $Y$ is an elliptic curve and the result is known; we will thus assume $g\ge 2$.
 
 Let $\{\pi^j_Y\}$ be the CK decomposition defined in \eqref{ck}. We observe that this CK decomposition is {\em generically defined\/} with respect to the family $\YY\to B$ (Notation \ref{not}), i.e. it is obtained by restriction from
 ``universal projectors'' $\pi^j_\YY\in A^{2g-1}(\YY\times_B \YY)$. This is just because $h$ and $\Delta_Y$ are generically defined.
 
  Writing $h^j(Y):=(Y,\pi^j_Y)\in\MM_{\rm rat}$, we have
   \begin{equation}\label{h2j} h^{2j}(Y) \cong \one(-j)\ \ \  \hbox{in}\ \MM_{\rm rat}\ \ \ \ \ \ (j=0,\ldots, 2g-1)\ .\end{equation}

 What we need to prove is that this CK decomposition is MCK, i.e.
      \begin{equation}\label{this} \pi_Y^k\circ \Delta_Y^{sm}\circ (\pi_Y^i\times \pi_Y^j)=0\ \ \ \hbox{in}\ A^{4g-2}(Y\times Y\times Y)\ \ \ \hbox{for\ all\ }i+j\not=k\ ,\end{equation}
      or equivalently that
       \[   h^i(Y)\otimes h^j(Y)\ \xrightarrow{\Delta_Y^{sm}}\  h(Y) \]
       coincides with 
       \[ h^i(Y)\otimes h^j(Y)\ \xrightarrow{\Delta_Y^{sm}}\ h(Y)\ \to\ h^{i+j}(Y)  \ \to\ h(Y)\ , \]   
       for all $i,j$.
              
   As a first step, let us assume that we have three integers $(i,j,k)$ with $i+j\not=k$ and at most one of them is odd. The cycle in \eqref{this} is generically defined and homologically trivial (cf. Remark \ref{rem:mck}). The isomorphisms \eqref{h2j} induce an injection
   \[ ({}^t \pi^i_Y\times{}^t \pi^j_Y\times\pi^k_Y)_\ast A^{4g-2}(Y\times Y\times Y)\ \hookrightarrow\ A^\ast(Y)\ ,\]
   and sends generically defined cycles to generically defined cycles (this is because the isomorphisms \eqref{h2j} are generically defined). As a consequence, the required vanishing \eqref{this} follows from the Franchetta property for $Y$, which is Proposition \ref{Fr}.      
   
   In the second step, let us assume that among the three integers $(i,j,k)$, exactly two are equal to $2g-1$. In this case, using the isomorphisms \eqref{h2j} we find an injection
   \[ ({}^t \pi^i_Y\times{}^t \pi^j_Y\times\pi^k_Y)_\ast A^{4g-2}(Y\times Y\times Y)\ \hookrightarrow\ A^\ast(Y\times Y)\ ,\] 
respecting the generically defined cycles. As such,  
the required vanishing \eqref{this} follows from the Franchetta property for $Y\times Y$, which is Proposition \ref{Fr}.
 
 In the third and final step we treat the case $i=j=k=2g-1$. In this case, we exploit the relation between $Y$ and the Fano variety $F:=F(Y)$. 
 Let $P\subset
 F\times Y$ denote the incidence correspondence (i.e. the universal $g-1$-dimensional linear subspace). As we have seen (Theorem \ref{YF}), $P\in A^{g}(F\times Y)$
  induces an isomorphism of motives
     \begin{equation}\label{chowiso}   {}^t P\colon\ \ h^{2g-1}(Y)\ \xrightarrow{\cong}\ h^1(F)(1-g)\ \ \ \hbox{in}\ \MM_{\rm rat}\ .\end{equation}
 The proof of this third and final step is split into two cases: the case where $g$ is at least $3$, and the case $g=2$.
     
    First, let us treat the case $g\ge 3$ (i.e. $\dim Y$ is at least $5$). In this case, $F$ being an abelian variety of dimension $g \ge 3$, it follows that there is a split injection
     \[  \sym^3 h^{2g-1}(Y)\cong \sym^3 h^1(F)(3-3g)\ \hookrightarrow\  h(F)(3-3g)\ \ \ \hbox{in}\ \MM_{\rm rat}\ .\]    
   Taking Chow groups, this induces a split injection
   \[  (\pi_Y^{2g-1}\times   \pi_Y^{2g-1}\times\pi_Y^{2g-1})_\ast A^{4g-2}(Y\times Y\times Y)^{\Sy_3}\ \hookrightarrow\ A^{g+1}(F) =0\ .\]
   (Here $\Sy_3$ is the symmetric group on 3 elements acting on $Y^3$ by permutation of the factors.) The small diagonal $\Delta_Y^{sm}$ is an element of $A^{4g-2}(Y\times Y\times Y)^{\Sy_3}$, and so
   \[   \pi_Y^{2g-1}\circ \Delta_Y^{sm}\circ (\pi_Y^{2g-1}\times \pi_Y^{2g-1})=    (\pi_Y^{2g-1}\times   \pi_Y^{2g-1}\times\pi_Y^{2g-1})_\ast (\Delta_Y^{sm})=   0\ \ \ \hbox{in}\ A^{4g-2}(Y\times Y\times Y)\   .\]  
     
 It remains to treat the case $g=2$ (i.e. $Y$ is a Fano threefold). In this case, $F$ is an abelian surface and so the isomorphism \eqref{chowiso} induces a split injection
   \[    \sym^2 h^{3}(Y)\cong \sym^2 h^1(F)(-2)\ \hookrightarrow\  h(F)(-2)\ \ \ \hbox{in}\ \MM_{\rm rat}\ .\]        
  Taking the tensor product with $h(Y)$, we get a split injection
    \[   \sym^2 h^3(Y)\otimes h(Y)\ \hookrightarrow\  h(F)(-2)\otimes h(Y)\ \ \ \hbox{in}\ \MM_{\rm rat}\ .\]     
 This induces in particular a split injection of Chow groups
   \begin{equation}\label{injch} (\pi^3_Y\times\pi^3_Y\times\pi^3_Y)_\ast A^j(Y\times Y\times Y)^{\Sy_2}\ \hookrightarrow\ A^{j-2}(F\times Y)\ \end{equation}
   (where $\Sy_2$ is the symmetric group on 2 elements acting as permutation on the first 2 factors).
   Since the Abel--Jacobi isomorphism is generically defined, the injection \eqref{injch} respects generically defined cycles.
   The cycle
     \[  \pi_Y^3\circ \Delta_Y^{sm}\circ (\pi_Y^3\times \pi_Y^3)  = ( \pi^3_Y\times\pi^3_Y\times\pi^3_Y)_\ast (\Delta_Y^{sm})\ \ \in\    A^6(Y\times Y\times Y)^{\Sy_2} \]
     is generically defined, and homologically trivial. The injection \eqref{injch}, plus the Franchetta property for $F\times Y$ (Proposition \ref{Fr2}), then guarantees that this cycle is rationally trivial, i.e.
     \[    \pi_Y^3\circ \Delta_Y^{sm}\circ (\pi_Y^3\times \pi_Y^3)=0\ \ \ \hbox{in}\ A^6(Y\times Y\times Y)\ .\]
     The theorem is now proven.     
           \end{proof}

 \begin{remark} Looking at the above proof, one observes that the case $g>2$ is easier than the case $g=2$. The reason is that for $g>2$ the relation
 \[ \pi^{2g-1}_Y\circ \Delta^{sm}_Y\circ (\pi^{2g-1}_Y\times \pi^{2g-1}_Y)=0\ \ \ \hbox{in}\ A^{4g-2}(Y^3) \]
 is trivially fulfilled for reasons of dimension, whereas in case $g=2$ one needs some extra argument (involving $F\times Y$) to obtain this relation. 
 \end{remark}

 \subsection{Compatibility} 
  In case $\dim Y=3$, the MCK decomposition for $Y$ behaves well with respect to the universal line on $Y$:
  
  \begin{proposition}\label{YandF} Let $Y\subset\PP^5$ be a smooth threefold that is a complete intersection of 2 quadrics. Let $F:=F(Y)$ be the Fano surface of lines on $Y$ and let $P\in A^2(F\times Y)$ be the universal line.
Then
  \[ P\ \ \in\ A^2_{(0)}(F\times Y)\ ,  \]
  where $F\times Y$ is given the product MCK decomposition. 
  In particular, 
    \[ P_\ast A^\ast_{(j)}(F)\ \subset\ A^\ast_{(j)}(Y)\ \ \ \hbox{and}\ \ P^\ast A^\ast_{(j)}(Y)\ \subset\ A^\ast_{(j)}(F)\  .\]
    \end{proposition}
    
    \begin{proof} The universal line $P$, and the MCK decomposition for $Y$ are generically defined (with respect to the base $B$). The MCK decomposition for the abelian surface $F$ is also generically defined with respect to $B$ (this is because the CK decomposition of \cite{DM}, which is generically defined with respect to any given family of abelian varieties, is MCK, cf. \cite{SV}). The result now follows from the Franchetta property for $F\times Y$ (Proposition \ref{Fr2}).
  \end{proof}
  
  \begin{remark} It seems likely that Proposition \ref{YandF} is true for $Y\subset\PP^{2g+1}$ a complete intersection of 2 quadrics with $g$ arbitrary. I have not been able to prove this.
   \end{remark}

 \section{The tautological ring}
 
 \begin{corollary}\label{cor1} Let $Y\subset\PP^{2g+1}$ be a smooth dimensionally transverse intersection of 2 quadrics, and let $m\in\NN$. Let
  \[ R^\ast(Y^m):=\Bigl\langle (p_i)^\ast(h), (p_{ij})^\ast(\Delta_Y)\Bigr\rangle\ \subset\ \ \ A^\ast(Y^m)   \]
  be the $\QQ$-subalgebra generated by pullbacks of the polarization $h\in A^1(Y)$ and pullbacks of the diagonal $\Delta_Y\in A^{2g-1}(Y\times Y)$. (Here $p_i$ and $p_{ij}$ denote the various projections from $Y^m$ to $Y$ resp. to $Y\times Y$).
  The cycle class map induces injections
   \[ R^\ast(Y^m)\ \hookrightarrow\ H^\ast(Y^m,\QQ)\ \ \ \hbox{for\ all\ }m\in\NN\ .\]
   \end{corollary}

\begin{proof} This is inspired by the analogous result for cubic hypersurfaces \cite[Section 2.3]{FLV3}, which in turn is inspired by analogous results for hyperelliptic curves \cite{Ta2}, \cite{Ta} (cf. Remark \ref{tava} below) and for K3 surfaces \cite{Yin}.

As in \cite[Section 2.3]{FLV3}, let us write $o:={1\over 4} h^{2g-1}\in A^{2g-1}(Y)$, and
  \[ \tau:= \Delta_Y - {1\over 4}\, \sum_{j=0}^{2g-1}  h^j\times h^{2g-1-j}\ \ \in\ A^{2g-1}(Y\times Y) \]
  (this cycle $\tau$ is nothing but the projector on the motive $h^{2g-1}(Y)$ considered above).
Moreover, let us write 
  \[ \begin{split}   o_i&:= (p_i)^\ast(o)\ \ \in\ A^{2g-1}(Y^m)\ ,\\
                        h_i&:=(p_i)^\ast(h)\ \ \in \ A^1(Y^m)\ ,\\
                         \tau_{ij}&:=(p_{ij})^\ast(\tau)\ \ \in\ A^{2g-1}(Y^m)\ .\\
                         \end{split}\]
We now define the $\QQ$-subalgebra
  \[ \bar{R}^\ast(Y^m):=\Bigl\langle o_i, h_i, \tau_{ij}\Bigr\rangle\ \ \ \subset\ H^\ast(Y^m,\QQ)\ \ \ \ \ (1\le i\le m,\ 1\le i<j\le m)\ ; \]
 this is the image of $R^\ast(Y^m)$ in cohomology. One can prove (just as \cite[Lemma 2.12]{FLV3} and \cite[Lemma 2.3]{Yin}) that the $\QQ$-algebra $ \bar{R}^\ast(Y^m)$
  is isomorphic to the free graded $\QQ$-algebra generated by $o_i,h_i,\tau_{ij}$, modulo the following relations:
    \begin{equation}\label{E:X'}
			o_i\cdot o_i = 0, \quad h_i \cdot o_i = 0,  \quad 
			h_i^{2g-1} =4\,o_i\,;
			\end{equation}
			\begin{equation}\label{E:X2'}
			\tau_{i,j} \cdot o_i = 0 ,\quad \tau_{i,j} \cdot h_i = 0, \quad \tau_{i,j} \cdot \tau_{i,j} = -b\, o_i\cdot o_j
			\,;
			\end{equation}
			\begin{equation}\label{E:X3'}
			\tau_{i,j} \cdot \tau_{i,k} = \tau_{j,k} \cdot o_i\,;
			\end{equation}
			\begin{equation}\label{E:X4'}
			\sum_{\sigma \in \mathfrak{S}_{b+2}} 
			\prod_{i=1}^{b/2+1} \tau_{\sigma(2i-1), \sigma(2i)} = 0\, ,
			\end{equation}
			where $b:= \dim H^{2g-1}(Y,\QQ)= 2g$.

To prove Corollary \ref{cor1}, it suffices to check that these relations are verified modulo rational equivalence.
The relations \eqref{E:X'} take place in $R^\ast(Y)$ and so they follow from the Franchetta property for $Y$ (Proposition \ref{Fr}). 
The relations \eqref{E:X2'} take place in $R^\ast(Y^2)$. The first and the last relations are trivially verified, because ($Y$ being Fano) one has
$A^{4g-2}(Y^2)=\QQ$. As for the second relation of \eqref{E:X2'}, this is relation \eqref{hyp}.
   
   Relation \eqref{E:X3'} takes place in $R^\ast(Y^3)$ and follows from the MCK decomposition (Theorem \ref{main}). Indeed, we have
   \[  \Delta_Y^{sm}\circ (\pi^{2g-1}_Y\times\pi^{2g-1}_Y)=   \pi^{4g-2}_Y\circ \Delta_Y^{sm}\circ (\pi^{2g-1}_Y\times\pi^{2g-1}_Y)  \ \ \ \hbox{in}\ A^{4g-2}(Y^3)\ ,\]
   which (using Lieberman's lemma) translates into
   \[ (\pi^{2g-1}_Y\times \pi^{2g-1}_Y\times\Delta_Y)_\ast    \Delta_Y^{sm}  =   ( \pi^{2g-1}_Y\times \pi^{2g-1}_Y\times\pi^{4g-2}_Y)_\ast \Delta_Y^{sm}                                                         \ \ \ \hbox{in}\ A^{4g-2}(Y^3)\ ,\]
   which means that
   \[  \tau_{13}\cdot \tau_{23}= \tau_{12}\cdot o_3\ \ \ \hbox{in}\ A^{4g-2}(Y^3)\ .\]
   
  Finally, relation \eqref{E:X4'}, which takes place in $R^\ast(Y^{b+2})$
   follows from the Kimura finite dimensionality relation \cite{Kim}: relation \eqref{E:X4'} expresses the vanishing
    \[ \sym^{b+2} H^{2g-1}(Y,\QQ)=0\ ,\]
    where $H^{2g-1}(Y,\QQ)$ is seen as a super vector space.
 This relation is also verified modulo rational equivalence, 
 (i.e., relation \eqref{E:X4'} is also true in $A^{(b+2)(2g-1)}(Y^{b+2})$): relation \eqref{E:X4'} involves a cycle in
   \[ A^\ast(\sym^{b+2} h^n(Y))\ ,\]
  and $\sym^{b+2} h^n(Y)$ is $0$ 
 because $Y$ has finite-dimensional motive (Lemma \ref{fd}). 
 This ends the proof.
%
%
%
 \end{proof}

\begin{remark}\label{tava} Given any curve $C$ and an integer $m\in\NN$, one can define the {\em tautological ring\/}
  \[ R^\ast(C^m):=  \Bigl\langle  (p_i)^\ast(K_C),(p_{ij})^\ast(\Delta_C)\Bigr\rangle\ \ \ \subset\ A^\ast(C^m) \]
  (where $p_i, p_{ij}$ denote the various projections from $C^m$ to $C$ resp. $C\times C$).
  Tavakol has proven \cite[Corollary 6.4]{Ta} that if $C$ is a hyperelliptic curve, the cycle class map induces injections
    \[  R^\ast(C^m)\ \hookrightarrow\ H^\ast(C^m,\QQ)\ \ \ \hbox{for\ all\ }m\in\NN\ .\]
   On the other hand, there are many (non hyperelliptic) curves for which the tautological rings $R^\ast(C^2)$ or $R^\ast(C^3)$ do {\em not\/} inject into cohomology (this is related to the non-vanishing of the Faber--Pandharipande cycle of Remark \ref{fp}; it is also
   related to the non-vanishing of the Ceresa cycle, cf. \cite[Remark 4.2]{Ta} and \cite[Example 2.3 and Remark 2.4]{FLV2}). 
   
Corollary \ref{cor1} shows that the tautological ring of intersections of 2 quadrics behaves similarly to that of hyperelliptic curves. 
\end{remark}

%
%

 \vskip0.5cm
\begin{nonumberingt} Thanks to the referee for helpful comments. Thanks to Kai and Len and Yasuyo for making the 2020 Lock-Down a pleasurable experience.

\end{nonumberingt}

\end{document}